\documentclass[11pt]{article}

\usepackage{amsthm}
\usepackage{latexsym}
\usepackage{dsfont}
\usepackage{amsmath}
\usepackage{color}
\usepackage[noblocks,affil-it]{authblk}
\usepackage{bbm}
\usepackage{amssymb}
\usepackage{amsmath}
\usepackage{mathtools}
\usepackage{wrapfig}
\usepackage{caption}
\usepackage{subcaption}

\numberwithin{equation}{section}
\theoremstyle{plain}
\newtheorem{Thm}{Theorem}[section]
\newtheorem{Lemma}[Thm]{Lemma}

\newtheorem{Cor}[Thm]{Corollary}
\theoremstyle{definition}
\newtheorem{Exa}[Thm]{Example}

\newtheorem{Rem}[Thm]{Remark}

\newcommand{\defeq}{\vcentcolon=}
\newcommand{\eqdef}{=\vcentcolon}
\newcommand{\distto}{\stackrel{\mathrm{d}}{\to}}

\newcommand{\N}{\mathbb{N}}

\newcommand{\R}{\mathbb{R}}

\newcommand{\imag}{\mathrm{i}}
\newcommand{\Real}{\mathrm{Re}}

\newcommand{\Prob}{\mathbb{P}}
\DeclareMathOperator{\Var}{\mathbb{V}\mathrm{ar}}

\newcommand{\E}{\mathbb{E}}

\renewcommand{\L}{\mathcal{L}}

\newcommand{\1}{\mathbbm{1}}

\newcommand{\dt}{\mathrm{d} \mathit{t}}
\newcommand{\du}{\mathrm{d} \mathit{u}}
\newcommand{\dx}{\mathrm{d} \mathit{x}}

\begin{document}

\title{Moment convergence of first-passage times in renewal theory}

\author[1]{Alexander Iksanov\thanks{E-mail: iksan@univ.kiev.ua}}
\affil[1]{Faculty of Cybernetics, Taras Shevchenko National University of Kyiv, 01601 Kyiv, Ukraine}
\affil[2]{Institut f\"{u}r Mathematische Statistik, Westf\"{a}lische Wilhelms-Universit\"{a}t M\"{u}nster, 48149 M\"{u}nster, Germany}
\affil[3]{Department of Mathematics, Technical University of Darmstadt, 64289 Darmstadt, Germany}
\author[1,2]{Alexander Marynych\thanks{E-mail: marynych@unicyb.kiev.ua}}
\author[3]{Matthias Meiners\thanks{E-mail: meiners@mathematik.tu-darmstadt.de}}

\maketitle

\begin{abstract}
Let $\xi_1, \xi_2, \ldots$ be independent copies of a positive random variable $\xi$,
$S_0 = 0$, and $S_k = \xi_1+\ldots+\xi_k$, $k \in \N_0$.
Define $N(t) = \inf\{k \in \N: S_k>t\}$ for $t\geq 0$.
The process $(N(t))_{t\geq 0}$ is the first-passage time process associated with $(S_k)_{k\geq 0}$.
It is known that if the law of $\xi$ belongs to the domain of attraction of a stable law or $\Prob(\xi>t)$ varies slowly at $\infty$,
then $N(t)$, suitably shifted and scaled,
converges in distribution as $t \to \infty$ to a random variable $W$ with a stable law or a Mittag-Leffler law.
We investigate whether there is convergence of the power and exponential moments to the corresponding moments of $W$.
Further, the analogous problem for first-passage times of subordinators is considered.
\end{abstract}

{\it Keywords:} exponential moment; L\'{e}vy process; power moment; renewal process; subordinator

{\it MSC 2010:} Primary 60K05, Secondary 60F05

\section{Introduction and results}

\paragraph{Setup.}
Let $\xi_1, \xi_2, \ldots$ be independent copies of a
positive random variable $\xi$.
We set $\mu \defeq \E[\xi] \in (0,\infty]$, and then $\sigma^2 \defeq \Var[\xi]$ whenever $\mu$ is finite.
Throughout the paper, we assume that the law of $\xi$ is non-degenerate, that is, $\Prob(\xi = c) < 1$ for all $c > 0$.
Define
\begin{equation*}
S_0 \defeq  0,  \quad   S_k \defeq  \xi_1+\ldots+\xi_k, \quad   k \in \N,
\end{equation*}
and
\begin{equation*}
N(t)    \defeq    \#\{k\in\N_0: S_k \leq t\}  = \inf\{k\in\N: S_k > t\},    \quad   t \geq 0.
\end{equation*}
The stochastic process $(N(t))_{t \geq 0}$ is called \emph{first-passage time process} associated with $(S_k)_{k \geq 0}$. The term `renewal counting process' is also used.

\paragraph{Objective.}
It is known (see, for instance, \cite[Proposition A.1]{Gnedin_al:2009}) % \cite[Theorem 3b]{Bingham:1972} or \cite{Whitt:2002})
that if the law of $\xi$ is in the domain of attraction of a stable law or $\Prob(\xi>t)$ varies slowly at $\infty$,
then
\begin{equation}    \label{eq:distributional convergence of N}
\frac{N(t)-b(t)}{a(t)} ~\distto~    W	\quad	\text{as}\quad t \to \infty
\end{equation}
where ``$\distto$'' denotes convergence in distribution,
$W$ is a non-degenerate random variable,
and $b(t) \in \R$, $a(t) > 0$ are suitable
shifting and scaling functions, respectively.

The purpose of this note is to answer the question: when does \eqref{eq:distributional convergence of N} imply
convergence of the corresponding power and exponential moments, finite or infinite? The motivation for writing
a short note on this problem comes from the fact that the moment convergence of first-passage time processes
repeatedly turned out to be an important technical step in other works on processes bearing some regenerative or
renewal structure. For instance, Theorems \ref{Thm:moment convergence A1,A2} and
\ref{Thm:moment convergence A3} below
are essential ingredients in our work on the finite-dimensional convergence of shot noise processes \cite{Iksanov+Marynych+Meiners:2014}.
Theorem \ref{Thm:moment convergence A4} is used to prove convergence of shot noise processes
to fractionally integrated inverse stable subordinators \cite{Iksanov+al:2016}.
Corollary \ref{Cor:moment convergence subordinator} is used in the proof of Theorem 3.3 in \cite{Gnedin+Iksanov:2012}.
Consequently, we found it useful to have one paper which contains the complete results on convergence
of power and exponential moments for renewal counting processes.

Before we state our results we briefly recall the different regimes
in which \eqref{eq:distributional convergence of N} holds.

\paragraph{Domains of attraction.}
The law of a random variable $\xi$ is in the domain of attraction of an $\alpha$-stable law, $\alpha\in (0,2]$ or $\Prob\{\xi>t\}$ varies slowly at $\infty$ if one of the following alternatives prevails\footnote{
Here, we do not treat the case where $\Prob(\xi>t)$ is regularly varying of index $-1$ at $\infty$
as it appears less frequently in applications and requires cumbersome calculations that would impair
the character of this paper as a brief note.
}:
\begin{itemize}\itemsep-2pt
	\item[(A1)]	$\mu < \infty$ and $\sigma^2 \defeq \Var[\xi] < \infty$;
	\item[(A2)]	$\mu < \infty$ but $\sigma^2 = \infty$ and $\ell_2(t) \defeq \E[\xi^2 \1_{\{\xi \leq t\}}]$ is slowly varying at $\infty$;
	\item[(A3)]	$\Prob(\xi > t) = t^{-\alpha} \ell(t)$ for some $\alpha \in (1,2)$ and a function $\ell$ slowly varying at $\infty$;
	\item[(A4)] $\Prob(\xi > t) = t^{-\alpha} \ell(t)$ for some $\alpha \in [0,1)$ and a function $\ell$ slowly varying at $\infty$.
\end{itemize}
We refer to \cite[Section 2.6]{Ibragimov+Linnik:1971} for details. The convergence of the first-passage time process in \eqref{eq:distributional convergence of N}
can now be described more precisely:
\begin{itemize}\itemsep-2pt
	\item[(N1)]
		if (A1) holds, then $b(t) = t/\mu$, $a(t) = \sigma\mu^{-3/2} c(t)$, $c(t) = \sqrt{t}$,
		and $W$ is a standard normal random variable;
    \item[(N2)]
    		if (A2) holds, then $b(t) = t/\mu$, $a(t) = \mu^{-3/2}c(t)$
		where $c(t)$ is a positive function satisfying $\lim_{t \to \infty} t \ell_2(c(t))c(t)^{-2} = 1$,
		and $W$ is a standard normal random variable;
    \item[(N3)]
    		if (A3) holds, then $b(t) = t/\mu$, $a(t) = \mu^{-(1+\alpha)/\alpha} c(t)$
		where $c(t)$ is a positive function such that
        	$\lim_{t \to \infty} t \ell(c(t))c(t)^{-\alpha} = 1$,
		and $W$ is a random variable with characteristic function given by\footnote{For $\alpha\in (1,2)$, $\Gamma(1-\alpha)$ is understood as $-\Gamma(2-\alpha)/(\alpha-1)$.}
		 \begin{equation}    \label{eq:characteristic function}
		\psi(\lambda)
		=  \exp\big\{-|\lambda|^\alpha \Gamma(1-\alpha)(\cos(\pi\alpha/2)+ \imag \sin(\pi\alpha/2)\, {\rm sgn}(\lambda))\big\},\ \lambda \in \R
		\end{equation}
		where $\Gamma(\cdot)$ denotes Euler's gamma function;
	\item[(N4)]
		if (A4) holds, then $b(t)=0$, $a(t) = 1/\Prob(\xi>t)$,
		and $W$ has a Mittag-Leffler distribution with parameter $\alpha$ (exponential with mean $1$ if $\alpha=0$), that is,
		$W$ has moment generating function
		\begin{center}	\label{eq:Mittag-Leffler law}
		$\E[e^{\theta W}]
		= E_{\alpha}\big(\frac{\theta}{\Gamma(1-\alpha)}\big)<\infty,
		\quad	\theta \in \R$
		\end{center}
		where here and throughout the paper,
		 $E_{\alpha}$ is the Mittag-Leffler function with parameter $\alpha$
		 given by $E_\alpha(z) \defeq \sum_{k\geq 0}\frac{z^k}{\Gamma(k\alpha+1)}$ for $z \in \R$.
\end{itemize}
%We refer to \cite{Gut:2009,Whitt:2002} for the cases (A1)-(A3), and to \cite[Corollary 3.4]{Meerschaert+Scheffler:2004} for (A4).

\paragraph{Main results for random walks.}
In what follows we use the notation $x_{-}$ and $x_{+}$ for the negative and positive part of a real number $x$:
\begin{equation*}
x_{-}:=-\min\{x,0\}\quad\text{and}\quad x_{+}:=\max\{x,0\}.
\end{equation*}

\begin{Thm}	\label{Thm:moment convergence A1,A2}
Suppose that either (A1) or (A2) holds, i.\,e., $\mu < \infty$ and either $\sigma^2 < \infty$
or $\sigma^2 = \infty$
and $\ell_2(t) \defeq \E[\xi^2 \1_{\{\xi \leq t\}}]$ is slowly varying at $\infty$.
Then
\begin{equation}	\label{eq:exponential moment convergence A1}
\lim_{t \to \infty} \E \bigg[\exp\bigg(\theta \frac{N(t)-t/\mu}{a(t)}\bigg)\bigg]
= \E[e^{\theta W}]
= e^{\frac{\theta^2}{2}},
\quad	\text{for every } \theta \geq 0
\end{equation}
where $W$ is standard normal,
$a(t) = \sqrt{\sigma^2\mu^{-3}t}$ in the case (A1) and $a(t) = \mu^{-3/2}c(t)$
for a positive function $c(t)$ satisfying $\lim_{t \to \infty} t \ell_2(c(t))c(t)^{-2} = 1$
in the case (A2). In particular,	
\begin{equation}	\label{eq:power moment convergence A1}
\lim_{t \to \infty} \E \bigg[\bigg(\frac{N(t)-t/\mu}{a(t)}\bigg)_{+}^p\bigg]
= \E [W_{+}^p]
= \frac{2^{p/2-1} \Gamma(\frac{p+1}{2})}{\sqrt{\pi}}
\quad	\text{for every } p>0.
\end{equation}
Further, in the case (A1)
\begin{equation}	\label{eq:power negative moment convergence A1}
\lim_{t \to \infty} \E \bigg[\bigg(\frac{N(t)-t/\mu}{a(t)}\bigg)_{-}^p\bigg]
= \E [W_{-}^p]
= \frac{2^{p/2-1} \Gamma(\frac{p+1}{2})}{\sqrt{\pi}},
\end{equation}
for every  $p\in [0,2]$. In the case (A2) the relation \eqref{eq:power negative moment convergence A1} holds for $p\in [0,2)$ and
\begin{equation}	\label{eq:power negative moment convergence A2}
\E[(N(t)-t/\mu)^2] ~\sim ~ \frac{2t}{\mu^3}\int_0^{t}\left(\int_x^{\infty}\Prob(\xi>z){\rm d}z\right) \dx,
\quad	\text{as } t \to \infty.
\end{equation}
\end{Thm}

\begin{Rem}	\label{Rem:moment convergence A1}
Without further assumptions on the law $\xi$,
the result stated in Theorem \ref{Thm:moment convergence A1,A2} is best possible in the following sense.
There exists a law for $\xi$ such that $\E[\xi^2]<\infty$ and
\begin{equation}	\label{eq:power negative moment no convergence A1}
\lim_{t \to \infty} \E \bigg[\bigg(\frac{N(t)-t/\mu}{\sqrt{t}}\bigg)_{-}^p\bigg]=\infty,
\end{equation}
for every $p>2$. An example is provided at the end of Section \ref{subsec:Proofs A1,A2,A3}.
\end{Rem}

\begin{Rem}
Convergence \eqref{eq:power negative moment convergence A1} is well-known in the case (A1) (see, for instance, \cite[Theorem 3.8.4]{Gut:2009}).
The asymptotic relation \eqref{eq:power negative moment convergence A2} follows from \cite[Theorems 2.3 and 2.4]{Mohan:1976}.
\end{Rem}

\begin{Thm}	\label{Thm:moment convergence A3}
Suppose that (A3) holds, i.\,e., $\Prob(\xi > t) = t^{-\alpha} \ell(t)$ for some $\alpha \in (1,2)$ and some $\ell$ slowly varying at $\infty$.
Let $c(t)$ be a positive function such that $\lim_{t \to \infty} t \ell(c(t))c(t)^{-\alpha} = 1$
and let $W$ be a random variable with characteristic function given by \eqref{eq:characteristic function}.
Then, for every $\theta \geq 0$, we have
\begin{equation}	\label{eq:exponential moments A3}
\lim_{t\to\infty} \E\bigg[\exp\bigg(\theta \frac{N(t)-t/\mu}{a(t)}\bigg)\bigg]
= \E[e^{\theta W}] = e^{-\Gamma(1-\alpha) \theta^{\alpha}},
\end{equation}
where $a(t) = \mu^{-(1+\alpha)/\alpha} c(t)$.
Further,
\begin{equation}	\label{eq:power moments A3}
\lim_{t\to\infty} \frac{\E[(N(t)-t/\mu)_{\pm}^p]}{a(t)^p}
= \E [W_{\pm}^{p}]
\quad	\text{for all } p > 0,
\end{equation}
where $\E[W_+^p] < \infty$ for all $p > 0$ and $\E[W_-^p] < \infty$ if and only if $p<\alpha$.
In particular,
\begin{align*}
\lim_{t\to\infty} & \frac{\E[|N(t)-t/\mu|^{p}]}{a(t)^p}
=\E [|W|^{p}]	\\
&=
\begin{cases}
\frac{2\Gamma(p+1)}{\pi p} \sin \big(\frac{\pi p}{2}\big) \Gamma\big(1\!-\!\frac p\alpha\big)
|\Gamma(1\!-\!\alpha)|^{\frac p\alpha} \cos\big(\frac{\pi p}{2}\!-\!\frac{\pi p}{\alpha}\big)	&	\text{for } 0 < p < \alpha,\\
\infty	&	\text{for } p \geq \alpha.
\end{cases}
\end{align*}
\end{Thm}

\begin{Thm}	\label{Thm:moment convergence A4}
Suppose that (A4) holds, i.\,e.,
$\Prob(\xi > t) = t^{-\alpha} \ell(t)$ for some $\alpha \in [0,1)$ and some $\ell$ slowly varying at $\infty$.
Then
\begin{equation}	\label{eq:exponential moment convergence A4}
\lim_{t\to\infty} \E\big[e^{\theta \Prob(\xi>t) N(t)}\big] = \E\big[e^{\theta W}\big]
= E_{\alpha}\Big(\frac{\theta}{\Gamma(1-\alpha)}\Big)<\infty
\text{ for every } \theta \in \R.
\end{equation}
In particular,
\begin{equation}	\label{eq:power moment convergence A4}
\lim_{t\to\infty} \E\big[(\Prob(\xi>t) N(t))^p\big] = \E\big[W^p\big]
= \frac{\Gamma(p+1)}{\Gamma(1+\alpha)^p \Gamma(p\alpha+1)} < \infty
\text{ for every } p \geq 0.
\end{equation}
\end{Thm}

\paragraph{Main results for subordinators.} Let $(X_t)_{t \geq 0}$ denote a subordinator, i.\,e., a nondecreasing L{\'e}vy process,
with $X_0=0$, drift coefficient $m \geq 0$, no killing and L{\'e}vy measure $\Pi$
that is concentrated on $\R^{+} \defeq [0,\infty)$.
Notice that compound Poisson processes are not excluded.
Put,
\begin{center}
$T_r \defeq \inf\{t \geq 0: X_t>r\},	\quad	r \geq 0.$
\end{center}
The stochastic process $(T_r)_{r \geq 0}$ is called \emph{first-passage time process} associated with $(X_t)_{t \geq 0}$.
The counterpart of \eqref{eq:distributional convergence of N} for $(T_r)_{r\geq 0}$ is % the first-passage time process is
\begin{equation}    \label{eq:distributional convergence of T_r}
\frac{T_r-b(r)}{a(r)} ~\stackrel{\mathrm{d}}{\to}~ W	\quad	\text{as}\quad 	r \to \infty
\end{equation}
for suitable constants $b(r) \in \R$ and $a(r) > 0$.
Let $N_r \defeq \inf \{k \in \N: X_k> r\}$ for $r \geq 0$.
Then $(N_r)_{r \geq 0}$ is the first-passage time process of $(X_n)_{n \in \N_0}$.
Clearly, $T_r \leq N_r \leq T_r+1$.
Hence, \eqref{eq:distributional convergence of T_r} holds if and only if
\eqref{eq:distributional convergence of N} holds with $N(t)$ replaced by $N_t$.
Furthermore, convergence of exponential or power moments in \eqref{eq:distributional convergence of T_r} holds if, and only if,
the corresponding convergence for the moments of $N_r$ holds.
We summarize these observations in the following corollary.

\begin{Cor}	\label{Cor:moment convergence subordinator}
Let $(X_t)_{t \geq 0}$ be a subordinator with $X_0=0$,
drift coefficient $m \geq 0$, no killing and L{\'e}vy measure $\Pi$ concentrated on $\R^{+}$.
Define $\xi \defeq X_1$, $\mu \defeq \E[\xi]$ and $\sigma^2 \defeq \Var[\xi]$. Then the following assertions hold:
\begin{itemize}
	\item[(a)]
		If the law of $\xi$ satisfies (A1), equivalently, $\int_{\{|x| \geq 1\}} x^2 \, \Pi(\dx) < \infty$, then
		\begin{equation}	\label{eq:exponential moment convergence A1 subordinator}
		\lim_{r \to \infty} \E \bigg[\exp\bigg(\theta \frac{T_r-r/\mu}{a(r)}\bigg)\bigg]
		= \E[e^{\theta W}]
		= e^{\frac{\theta^2}{2}}
		\quad	\text{for every } \theta\geq 0
		\end{equation}
		where $W$ is standard normal and $a(r):=\sqrt{\sigma^2\mu^{-3}r}$. In particular,
		\begin{equation}	\label{eq:power moment convergence A1 subordiantor}
\lim_{r \to \infty} \E \bigg[\bigg(\frac{T_r-r/\mu}{a(r)}\bigg)_{+}^p\bigg] = \E [W_{+}^p]
= \frac{2^{p/2-1} \Gamma(\frac{p+1}{2})}{\sqrt{\pi}}
\quad	\text{for every } p>0.
\end{equation}
Further,
\begin{equation}	\label{eq:power negative moment convergence A1 subordinator}
\lim_{r \to \infty} \E \bigg[\bigg(\frac{T_r-r/\mu}{a(r)}\bigg)_{-}^p\bigg]
= \E [W_{-}^p]
= \frac{2^{p/2-1} \Gamma(\frac{p+1}{2})}{\sqrt{\pi}}
\quad	\text{for every } p\in [0,2].
\end{equation}
		
	\item[(b)]
		If the law of $\xi$ satisfies (A2), equivalently,
		\begin{equation}\label{eq:A2'}
		\mu<\infty,\;\sigma^2 = \infty
		\text{ and }
		\ell_2^\Pi(t) \defeq \! \int_{(1,t]} \!\!\!\!\! x^2 \, \Pi(\dx)\;\text{is slowly varying at }\infty,
		\end{equation}
		then \eqref{eq:exponential moment convergence A1 subordinator}, \eqref{eq:power moment convergence A1 subordiantor} and \eqref{eq:power negative moment convergence A1 subordinator} (the latter only for $0<p<2$) hold with $a(r) = \mu^{-3/2}c(r)$
		where $c(r)$ is a positive function satisfying $\lim_{r \to \infty} r \ell_2^\Pi(c(r))c(r)^{-2} = 1$.
	\item[(c)]
		If the law of $\xi$ satisfies (A3), equivalently, $\Pi((t,\infty)) \!=\! t^{-\alpha} \ell^\Pi(t)$ for some $\alpha \in (1,2)$
		and some $\ell^\Pi$ slowly varying at $\infty$,
		then \eqref{eq:exponential moment convergence A1 subordinator} holds with
		$a(r) = \mu^{-(1+\alpha)/\alpha}c(r)$
		where $c(r)$ is a positive function satisfying $\lim_{r \to \infty}  r\ell^\Pi(c(r))c(r)^{-\alpha}=  1$,
		and $W$ is a random variable with characteristic function given by \eqref{eq:characteristic function}.
		Further,
		\begin{equation}	\label{eq:power moments A3 subordinator}
		\lim_{r\to\infty} \frac{\E[(T_r-r/\mu)_{\pm}^p]}{a(r)^p}
		= \E [W_{\pm}^{p}]
		\quad	\text{for all } p > 0.
		\end{equation}
	\item[(d)]
		If $\xi$ satisfies (A4), equivalently, $\Pi((t,\infty)) = t^{-\alpha} \ell^\Pi(t)$ for some $\alpha \in [0,1)$
		and a function $\ell^\Pi$ which is slowly varying at $\infty$, then
		\begin{equation}	\label{eq:exponential moment convergence A4 subordinator}
		\lim_{r\to\infty} \E\big[e^{\theta \Pi((r,\infty)) T_r}\big] = \E\big[e^{\theta W}\big]
		= E_{\alpha}\Big(\frac{\theta}{\Gamma(1-\alpha)}\Big)<\infty
		\text{ for every } \theta \in \R
		\end{equation}
		where $W$ has the Mittag-Leffler distribution with parameter $\alpha$ and $E_{\alpha}(\cdot)$ is the Mittag-Leffler function.
\end{itemize}
\end{Cor}

We close this section with a remark that the asymptotics of $\E[e^{aN(t)}]$ as $t\to\infty$ as well as exponential moments of the number of visits and the last-exit time was investigated in \cite{Iksanov+Meiners:2010b, Iksanov+Meiners:2010a} for random walks with two-sided jumps and in \cite{Aurzada+Iksanov+Meiners:2015} for L\'{e}vy processes.

%\paragraph{Literature.}
%Convergence of power moments in \eqref{eq:distributional convergence of N} has been studied before.
%We refer to \cite{Gut:2009} for an introduction and an overview over the literature.
%
%The asymptotic behavior as $t \to \infty$ of $\E[e^{a N(t)}]$ for $a>0$ has been investigated in \cite{Iksanov+Meiners:2010b}.
%However, in \cite{Iksanov+Meiners:2010b}, $a$ is not allowed to depend on $t$ and hence the results obtained there
%do not apply here.

\section{Proofs of the main results}

%\paragraph{Notation.}
We denote by $\varphi$ the Laplace transform of $\xi$, i.\,e., $\varphi(\lambda)=\E[e^{-\lambda \xi}]$, $\lambda \geq 0$.
Some relevant results about the behavior of $\varphi$ at $0$ which we use in the proofs below are collected in the Appendix.

%\noindent
%For positive sequences $(a_n)_{n \in \N}$ and $(b_n)_{n \in \N}$ we write
%\vspace{-0.2cm}
%\begin{itemize}\itemsep-1pt
%	\item
%		$a_n ~\sim ~ b_n$ as $n \to \infty$ if $\lim_{n \to \infty} a_n/b_n = 1$;
%	\item
%		$a_n = O(b_n)$ as $n \to \infty$ or $a_n \preccurlyeq b_n$ if $\limsup_{n \to \infty} \frac{a_n}{b_n} < \infty$;
%	\item
%		$a_n \asymp b_n$ as $n \to \infty$ if $a_n = O(b_n)$ and $b_n = O(a_n)$ as $n \to \infty$;
%	\item
%		$a_n = o(b_n)$ as $n \to \infty$ if $\lim_{n \to \infty} a_n/b_n = 0$.
%\end{itemize}

\subsection{Proof of Theorems \ref{Thm:moment convergence A1,A2}, \ref{Thm:moment convergence A3}}	\label{subsec:Proofs A1,A2,A3}

\paragraph{Convergence of exponential moments of positive order and power moments of the positive parts.}

In view of (N1), (N2) and (N3) we have
\begin{equation*}
e^{\theta \frac{N(t)-t/\mu}{a(t)}} \distto e^{\theta W}	\quad	\text{as } t \to \infty
\end{equation*}
for every $\theta\geq 0$, where $W$ is standard normal in the cases (A1) and (A2), and $W$ has characteristic function
given by \eqref{eq:characteristic function} in the case (A3).
Hence, it is enough to show that the family $(\exp(\theta a(t)^{-1}(N(t)-t/\mu)))_{t\geq t_0}$
is uniformly integrable for every $\theta>0$ and some $t_0>0$.
To this end, by the Vall\'{e}e-Poussin criterion of uniform integrability it suffices to check that
\begin{equation}	\label{eq:Vallee-Poussin}
\sup_{t \geq t_0} \E \Big[e^{\theta \frac{N(t)-t/\mu}{a(t)}}\Big]<\infty
\end{equation}
for every $\theta>0$.
While doing so, we can neglect the constant factors in the scaling functions $a(t)$ thus working with $c(t)$ in place of $a(t)$.
With the help of Markov's inequality we obtain
\begin{align}
\E \Big[e^{\theta \frac{N(t)-t/\mu}{c(t)}}\Big]
& =
\int_0^{\infty} \Prob\Big(e^{\theta \frac{N(t)-t/\mu}{c(t)}} > x\Big) \, \dx	\notag
=
\int_{-\infty}^{\infty} e^x \, \Prob\bigg(\theta \frac{N(t)-t/\mu}{c(t)} > x\bigg) \, \dx	\notag	\\
& \leq
1 + \int_{0}^{\infty} e^x \, \Prob\big(N(t) > x c(t)/\theta + t/\mu \big) \, \dx	\notag	\\
& =
1 + \int_{0}^{\infty} e^x \, \Prob\big(S_{\lfloor x c(t)/\theta + t/\mu \rfloor} \leq t\big) \, \dx	\notag	\\
& \leq
1+\int_{0}^{\infty} e^x \, e^{\lambda t} (\varphi(\lambda))^{x c(t)/\theta + t/\mu -1} \, \dx	\notag	\\
& \leq
1+(e^{\lambda\mu} \varphi(\lambda))^{t/\mu} \int_{0}^{\infty} e^x \, (\varphi(\lambda))^{x c(t) / \theta - 1} \, \dx
\label{eq:exponential moments integral estimate A1,A2}
\end{align}
for every $\lambda > 0$.
We will demonstrate that \eqref{eq:Vallee-Poussin} is a consequence of
\begin{equation}	\label{eq:Vallee-Poussin sufficient}
\sup_{t \geq t_0} \int_{0}^{\infty} e^x \, \varphi(\lambda/c(t))^{x c(t) / \theta - 1} \, \dx < \infty
\end{equation}
for some $\lambda>0$.

\noindent {\sc Case (A1)} in which $c(t) = \sqrt{t}$.
From formula \eqref{eq:varphi at 0 A1} in the Appendix we infer
\begin{equation*}
\varphi(\lambda/\sqrt{t}) = 1 - \frac{\mu\lambda}{\sqrt{t}} + \frac{\mu^2+\sigma^2}{2} \frac{\lambda^2}{t} + o\Big(\frac{1}{t}\Big)
\quad	\text{as } t \to \infty
\end{equation*}
whence
\begin{equation*}
e^{\lambda\mu/\sqrt{t}} \varphi(\lambda/\sqrt{t}) = 1+\frac{\sigma^2}{2} \frac{\lambda^2}{t} + o\Big(\frac{1}{t}\Big)
\quad	\text{as } t \to \infty.
\end{equation*}
Thus, substituting $\lambda$ by $\lambda/\sqrt{t}$ in \eqref{eq:exponential moments integral estimate A1,A2} we see that \eqref{eq:Vallee-Poussin sufficient} is indeed sufficient for \eqref{eq:Vallee-Poussin} to hold.
\smallskip

\noindent {\sc Case (A2)}.
From \eqref{eq:varphi at 0 A2} and the relation $\lim_{t \to \infty} t\ell_2(t)c(t)^{-2}=1$, we infer
\begin{equation*}
\varphi(\lambda/c(t))
= 1 - \frac{\mu\lambda}{c(t)} + \frac{1}{2} \frac{\lambda^2}{t} + o\Big(\frac1t\Big)
\quad	\text{as } t \to \infty
\end{equation*}
and, since $c(t)^{-2} = o(t^{-1})$ as $t \to \infty$,
\begin{equation*}
e^{\lambda\mu/c(t)} \varphi(\lambda/c(t)) = 1+\frac{1}{2} \frac{\lambda^2}{t} + o\Big(\frac{1}{t}\Big)
\quad	\text{as } t \to \infty.
\end{equation*}
Thus, substituting $\lambda$ by $\lambda/c(t)$ in \eqref{eq:exponential moments integral estimate A1,A2}
we conclude that \eqref{eq:Vallee-Poussin sufficient} is sufficient for \eqref{eq:Vallee-Poussin}.
\smallskip

\noindent {\sc Case (A3)}.
Since $c(t)\to\infty$ as $t\to\infty$, we infer from \eqref{eq:varphi at 0 A3} in the Appendix that
\begin{equation}	\label{eq:Taylor of varphi(lambda/c(t))}
\varphi(\lambda/c(t)) = 1-\frac{\mu\lambda}{c(t)} + \frac{\lambda^{\alpha}c_{\alpha}}{t} + o(1/t)
\end{equation}
whence
\begin{equation*}
e^{\lambda\mu/c(t)}\varphi(\lambda/c(t))=1+\frac{\lambda^{\alpha}c_{\alpha}}{t}+o(1/t)
\end{equation*}
as $t \to \infty$. This implies the claim.

It remains to prove \eqref{eq:Vallee-Poussin sufficient}. From formulae \eqref{eq:varphi at 0 A1}, \eqref{eq:varphi at 0 A2} and \eqref{eq:varphi at 0 A3} we deduce that for every fixed $\varepsilon \in (0,\mu)$, $\lambda > 0$ and sufficiently large $t$
\begin{equation*}
\varphi\Big(\frac{\lambda}{c(t)}\Big) \leq 1 - \frac{(\mu-\varepsilon)\lambda}{c(t)} \leq e^{-(\mu-\varepsilon)\lambda/c(t)}.
\end{equation*}
Consequently,
\begin{equation*}
\int_{0}^{\infty} e^x \, \varphi(\lambda/c(t))^{x c(t) / \theta - 1} \, \dx
\leq e^{(\mu-\varepsilon)\lambda/c(t)} \int_0^\infty e^{x(1-(\mu-\varepsilon)\lambda/ \theta)}\, \dx,
\end{equation*}
and the latter integral is finite provided that $\lambda$ is chosen large enough.

Thus, the first equalities in relations \eqref{eq:exponential moment convergence A1} and \eqref{eq:exponential moments A3} are proved.
The second equality in \eqref{eq:exponential moment convergence A1} is a well-known formula for exponential moments of a standard normal law.
The second equality in \eqref{eq:exponential moments A3}, namely, $\E[e^{\theta W}] = e^{-\Gamma(1-\alpha) \theta^{\alpha}}$ for all $\theta \geq 0$,
can be found in many sources, see, for instance, \cite[Exercise 29.15]{Sato:1999}.
Now the first parts of Theorems \ref{Thm:moment convergence A1,A2} and \ref{Thm:moment convergence A3} regarding the exponential moments of positive order are completely proved. Relations \eqref{eq:power moment convergence A1} and \eqref{eq:power moments A3} (the latter only for the positive parts) follow from the inequality $x^{p}_{+} \leq e^{p x}$ which yields the uniform integrability of the corresponding families.

\paragraph{Convergence of power moments of negative parts.}

We treat the cases (A2) and (A3) simultaneously. First fix $0<p<\alpha$ (with $\alpha=2$ in the case (A2)) and $r \in (p \vee 1,\alpha)$. As before it is enough to show that for some $t_0>0$,
\begin{equation*}
\sup_{t \geq t_0} \frac{\E[(N(t)-t/\mu)_{-}^{r}]}{c(t)^r} < \infty.
\end{equation*}
By the regular variation of $c$, this is implied by % equivalent to
\begin{equation}	\label{eq:sufficient negative power moments A3}
\E[(N(\mu n)-n)_{-}^{r}] = O(c(n)^r)	\quad	\text{as } n \to \infty.
\end{equation}
We have
\begin{align*}
\E&[(N(\mu n)-n)_{-}^{r}]	\\
&~=\sum_{k\geq 1} \Prob((N(\mu n)-n)_{-}^{r}\geq k)
= \sum_{k\geq 1} \Prob(N(\mu n)\leq  n-k^{1/r})	\\
&~= \sum_{k=1}^{\lfloor n^r \rfloor}\Prob(S_{\lfloor n-k^{1/r} \rfloor} > \mu n)
= \sum_{j=0}^{n-1}\sum_{k \in (j^r, (j+1)^r]} \Prob(S_{\lfloor n-k^{1/r} \rfloor} > \mu n)\\
&~\leq \sum_{j=0}^{n-1}((j+1)^{r}-j^{r})\Prob(S_{n-j-1} > \mu n)
\leq r\sum_{j=0}^{n-1}(j+1)^{r-1}\Prob(S_{n-j-1} > \mu n)\\
&~= r\sum_{j=1}^{n} j^{r-1}\Prob(S_{n-j} - (n-j)\mu > \mu j)
\leq r\sum_{j=1}^{n} j^{r-1}\Prob\Big(\max_{0 \leq i \leq n-1} (S_{i} - i\mu) > \mu j\Big)	\\
&~\leq  r + {\rm const} \cdot \E\Big[\max_{0 \leq i \leq n} |S_{i} - i\mu|^r\Big]
\leq  r + {\rm const} \cdot \E\big[|S_{n} - n\mu|^r\big]
= O(c(n)^r)
\end{align*}
as $n \to \infty$ where the penultimate step is a consequence of the maximal $\L^r$-inequality
and the last step follows from \cite[Lemma 5.2.2]{Ibragimov+Linnik:1971}. The formula for $\E [|W|^{p}]$, $0<p<\alpha$ in the case (A3) is justified by Lemma \ref{Lem:power moments A3}.

Finally, we show that in the case (A3) 
\begin{equation*}
\lim_{t \to \infty} \E\bigg[\frac{(N(t)-t/\mu)_{-}^p}{a(t)^p}\bigg] = \infty=\E[W_-^p]
\end{equation*}
for $p\geq \alpha$. The second equality follows from the well-known relation $\Prob\{W_->x\}\sim cx^{-\alpha}$ as $x\to\infty$ for a positive constant $c$. With this at hand the first equality is a consequence of (N3) and Fatou's lemma.
The proof of Theorems \ref{Thm:moment convergence A1,A2} and \ref{Thm:moment convergence A3} is complete.\smallskip

We close this section with an example showing that convergence of moments of order $p>2$
may fail in the case of a normal limit.

\begin{Exa}	\label{Exa:moment convergence A1}
If the survival function of $\xi$ is given by
\begin{equation*}
\Prob(\xi > t)=\frac{1}{(t+1)^2\log^2 (t+e)},\quad t\geq 0,
\end{equation*}
then $\E[\xi^2]<\infty$ and
\begin{equation*}
\Prob(S_n>\gamma n)\geq \Prob(\max\{\xi_1,\xi_2,\ldots,\xi_n\}>\gamma n)
=1-(1-\Prob(\xi>\gamma n))^n
\sim n\Prob(\xi >\gamma n),\quad n\to\infty
\end{equation*}
for every fixed $\gamma>0$. Therefore,
\begin{equation*}
\E [(N(2\mu n)-2n)_{-}^p] \geq n^p \Prob(N(2\mu n)\leq n)=n^p \Prob(S_n > 2\mu n)\geq c n^{p+1}\Prob(\xi > n)
\end{equation*}
for some $c>0$ and all sufficiently large $n$. Hence \eqref{eq:power negative moment no convergence A1} holds.
\end{Exa}

\paragraph{Alternative proof for the convergence of first absolute moments.}

There is an alternative elegant proof of the convergence of {\it the first moments} in \eqref{eq:distributional convergence of N}
for the cases (A1) through (A3) based on the representation
\begin{align*}
\E [|S_{N(\mu n)}-S_n|]
&=
\E[S_{N(\mu n)\vee n}-S_{N(\mu n) \wedge n}]   \\
&=
\mu\E[(N(\mu n)\vee n) - (N(\mu n)\wedge n)]
= \mu \E[|N(\mu n)-n|],
\end{align*}
where the second equality follows from Wald's identity.
From this one obtains
\begin{align}
\E[ |S_n-\mu n| - (S_{N(\mu n)} - \mu n)]
& =
\mu \E[|N(\mu n)-n|] \notag  \\
& \leq
\E[|S_n-\mu n| +(S_{N(\mu n)}-\mu n)]. \label{eq:bounds on E|N(mun)-n|}
\end{align}
According to \cite[Lemma 5.2.2]{Ibragimov+Linnik:1971}
\begin{equation}	\label{eq:RW convergence to stable 1st moment}
\lim_{n \to \infty} \frac{\E[|S_n-\mu n|]}{c(n)}  ~=~ \E[|W|].
\end{equation}
From \cite{Mohan:1976} it is known that, as $t \to \infty$,
\begin{equation*}
\E[S_{N(t)}-t]  \sim  \begin{cases}
                        \text{const}							&   \text{in the case (A1),}        \\
                        \text{const} \cdot \ell(t)				&   \text{in the case (A2),}        \\
                        \text{const} \cdot t^{2-\alpha} \ell(t)	&   \text{in the case (A3),}
                        \end{cases}
\end{equation*}
provided that the law of $\xi$ is non-lattice.

Assume now that the law of $\xi$ is lattice with span $d>0$.
In the case (A1), according to \cite[Theorem 9]{Feller:1949},
$\E[S_{N(nd)}-nd]$ tends to a constant as $n \to \infty$.
Hence
\begin{equation*}
\E[S_{N(t)}-t] = O(1) \quad \text{as } t \to \infty.
\end{equation*}
In the cases (A2) and (A3), according to \cite[Theorem 6]{Sgibnev:1981},
$\E[S_{N(t)}-t]$ exhibits the same asymptotic behavior as in the non-lattice case.

Since $c(t)$ is regularly varying of index $1/\alpha$ at $\infty$ (where $\alpha=2$ in the cases (A1) and (A2)),
we conclude that
\begin{equation*}
\lim_{n \to \infty} \frac{\E[S_{N(\mu n)}-\mu n]}{c(n)} = 0.
\end{equation*}
Applying this and \eqref{eq:RW convergence to stable 1st moment} to \eqref{eq:bounds on E|N(mun)-n|}
we infer
\begin{equation*}
\lim_{n \to \infty} \mu \frac{\E[|N(\mu n)-n|]}{c(n)}    ~=~ \E[|W|].
\end{equation*}
Now we have to check that this relation implies the convergence of the first absolute moments in \eqref{eq:distributional convergence of N}.
For any $t>0$ there exists an $n=n(t)\in\N_0$ such that $t \in (\mu n,\mu(n+1)]$.
Hence, by subadditivity,
\begin{equation*}
\E[N(t)-N(\mu n)] \leq \E[N(\mu(n+1))-N(\mu n)] \leq \E[N(\mu)].
\end{equation*}
It remains to observe that the regular variation of $c(t)$ entails $\lim_{t \to \infty} c(\mu n(t)\mu^{-1})/c(t)=\mu^{-1/\alpha}$.
This implies the asserted convergence of the first absolute moments in \eqref{eq:distributional convergence of N}.

\subsection{Proof of Theorem \ref{Thm:moment convergence A4}}

Arguing as in the proof of Theorems \ref{Thm:moment convergence A1,A2} and \ref{Thm:moment convergence A3},
we conclude that it suffices to show that
\begin{equation*}
\sup_{t \geq t_0} \E [e^{\theta \Prob(\xi>t) N(t)}]<\infty
\end{equation*}
for every $\theta>0$ and some $t_0\geq 0$. %For any fixed $t_0 > 0$, we have the bound $\E[e^{\theta \Prob(\xi>t)N(t)}] \leq \E[e^{\theta N(t_0)}] < \infty$
%by \cite[Theorem 2.1(b)]{Iksanov+Meiners:2010a}.
%Hence it suffices to show the finiteness of the supremum over $t \geq t_0$ for arbitrary large $t_0 > 0$. We infer
Write
\begin{align}
\frac{\E [e^{\theta \Prob(\xi>t)N(t)}]-1}{e^{\theta \Prob(\xi>t)}-1}
&= \sum_{k \geq 0}e^{\theta \Prob(\xi>t) k}\Prob(S_k \leq t)	\notag	\\
&\leq e^{\lambda t} \sum_{k \geq 0} e^{\theta \Prob(\xi>t) k} \varphi(\lambda)^k
=\frac{e^{\lambda t}}{1-e^{\theta \Prob(\xi>t)}\varphi(\lambda)}	\label{eq:sup bound A4}
\end{align}
for every $\lambda > 0$ such that $e^{\theta \Prob(\xi>t)} \varphi(\lambda) < 1$.
Pick an arbitrary
$c>(\theta/\Gamma(1-\alpha))^{1/\alpha}$
and note that
\begin{equation}	\label{eq:choice of lambda A4}
\frac{1-e^{-\theta \Prob(\xi>t)}}{1-\varphi(c/t)}
\sim \frac{\theta \Prob(\xi>t)}{\Gamma(1-\alpha)\Prob(\xi > t/c)}\to \frac{\theta
c^{-\alpha}}{\Gamma(1-\alpha)}<1
\end{equation}
as $t\to\infty$ where \eqref{eq:varphi at 0 A4} has been used.
Relation \eqref{eq:choice of lambda A4} entails \begin{equation*}e^{\theta \Prob(\xi>t)}\varphi(c/t)<1\end{equation*} for all $t>0$ large enough.
Therefore, choosing $\lambda = c/t$ in \eqref{eq:sup bound A4}
and using again \eqref{eq:choice of lambda A4} % using the asymptotic relation $1-\varphi(z) \sim \Gamma(1-\alpha)\Prob(\xi>z^{-1})$ as $z \downarrow 0$
we infer 
\begin{align*}
\E[e^{\theta \Prob(\xi>t)N(t)}] - 1
\leq e^c\frac{e^{\theta \Prob(\xi>t)}-1}{1-e^{\theta \Prob(\xi>t)}\varphi(c/t)}
\to \frac{e^c\theta}{\Gamma(1-\alpha)c^{\alpha}-\theta}
\quad \text{as } t \to \infty
\end{align*}
which completes the proof of the first equalities in \eqref{eq:exponential moment convergence A4} and \eqref{eq:power moment convergence A4}. While the second equality in
\eqref{eq:exponential moment convergence A4} and the second equality in \eqref{eq:power moment convergence A4} when $\alpha=0$ are immediate, the second equality in \eqref{eq:power moment convergence A4} when $\alpha\in (0,1)$ follows from Lemma \ref{Lem:power moments A4}. The proof of Theorem \ref{Thm:moment convergence A4} is complete.

\subsection{Proof of Corollary \ref{Cor:moment convergence subordinator}}
The claimed asymptotic relations follow almost immediately from Theorems \ref{Thm:moment convergence A1,A2}
to \ref{Thm:moment convergence A4} and the fact that $T_r \leq N_r \leq T_r+1$. It remains to check the claimed equivalent reformulations of (A1) through (A4) in terms of the L\'evy measure $\Pi$
and to make sure that we use the right scaling.\smallskip

\noindent
{\it Proof of (a):}
$\sigma^2 < \infty$ is equivalent to $\int_{\{|x| \geq 1\}} x^2 \, \Pi(\dx) < \infty$ by standard theory for L\'evy processes,
see \cite[Corollary 25.8]{Sato:1999}.	\smallskip

\noindent
{\it Proof of (b):}
In the proof of Lemma 6(a) in \cite[$3$rd line after (3.10)]{Embrechts+Goldie:1981},
it is shown that condition \eqref{eq:A2'} implies that
\begin{equation*}
\ell_2(t) = \E[\xi^2 \1_{\{\xi \leq t\}}] \sim \int_{(0,t]} x^2 \, \Pi(\dx) = \ell_2^\Pi(t)
\quad	\text{as }	t \to \infty.
\end{equation*}
Consequently, the asymptotic relation $\lim_{t \to \infty} t \ell_2(c(t))c(t)^{-2} = 1$ is equivalent to
\begin{equation*}
\lim_{r \to \infty} r\ell_2^\Pi(c(r))c(r)^{-2} = 1.
\end{equation*}\smallskip

\noindent
{\it (c) and (d):}
According to \cite[Proposition 0]{Embrechts+Goldie:1981},
$\Prob(\xi>t)$ is regularly varying of index $-\alpha$ at $\infty$ if and only if the same is true for
$\Pi((t,\infty))$, and in this case $\Prob(\xi>t) \sim \Pi((t,\infty))$ as $t \to \infty$.
This proves (d), while (c) follows upon noting that $\ell^\Pi(t) \sim \ell(t)$ which implies that
the asymptotic relations $\lim_{t \to \infty} t\ell(c(t))c(t)^{-\alpha} = 1$ and $\lim_{r \to \infty} r\ell^\Pi(c(r))c(r)^{-\alpha} = 1$
are equivalent.

\section{Appendix}

\subsection{Laplace transforms}
Here, we gather known results on the behavior of Laplace transforms at $0$ that play a role in our derivations.
Recall that $\varphi$ denotes the Laplace transform of $\xi$.
\smallskip

\noindent
In the case (A1), $\E[\xi^2] = \mu^2+\sigma^2 < \infty$ and a Taylor expansion of $\varphi$ at $0$ gives
\begin{equation}	\label{eq:varphi at 0 A1}
\varphi(\lambda) = 1 - \mu\lambda + \frac{\mu^2+\sigma^2}{2} \lambda^2 + o(\lambda^2)	\quad	\text{as } \lambda\to 0+.
\end{equation}

\noindent
In the case (A2), $\E[\xi^2] = \infty$ and $\ell_2(t) = \E[\xi^2 \1_{\{\xi \leq t\}}]$ is slowly varying at $\infty$.
Hence,
\begin{equation}	\label{eq:varphi at 0 A2}
\varphi(\lambda) - (1 - \mu \lambda)
\sim \frac{1}{2} \lambda^2 \ell_2(1/\lambda)
\quad	\text{as } \lambda \to 0+
\end{equation}
by the implication (8.1.11c)$\Rightarrow$ (8.1.9) of \cite[Theorem 8.1.6]{Bingham+Goldie+Teugels:1989}.

\noindent
In the case (A3), using that $\Prob(\xi > t)$ is regularly varying of index $-\alpha$ for $\alpha\in (1,2)$
we infer
\begin{equation}	\label{eq:varphi at 0 A3}
\varphi(\lambda)-(1-\mu \lambda) \sim c_{\alpha} \Prob(\xi > 1/\lambda)	\quad	\text{as } \lambda \to 0+
\end{equation}
with $c_{\alpha}\defeq \frac{\Gamma(2-\alpha)}{\alpha-1}$ by \cite[Theorem 8.1.6]{Bingham+Goldie+Teugels:1989}.\smallskip

\noindent
In the case (A4), since $\Prob(\xi > t)$ is regularly varying of index $-\alpha$ for $\alpha\in [0,1)$
an application of \cite[Corollary 8.1.7]{Bingham+Goldie+Teugels:1989} yields
\begin{equation}	\label{eq:varphi at 0 A4}
1-\varphi(\lambda) \sim \Gamma(1-\alpha)\Prob(\xi>1/\lambda)	\quad	\text{as } \lambda \to 0+.
\end{equation}

\subsection{Moment computations}

\begin{Lemma}		\label{Lem:power moments A3}
Let $W$ be a random variable with characteristic function given by \eqref{eq:characteristic function}.
Then, for $r < \alpha$,
\begin{equation*}
\E [|W|^r]
~=~ \frac{2\Gamma(r+1)}{\pi r} \sin \Big(\frac{\pi r}{2}\Big) \Gamma\Big(1-\frac r\alpha\Big)
|\Gamma(1-\alpha)|^{\frac r\alpha} \cos\Big(\frac{\pi r}{2}-\frac{\pi r}{\alpha}\Big).
\end{equation*}
In particular,
$\E [|W|]  = \frac2\pi \Gamma(1-\frac1\alpha)|\Gamma(1-\alpha)|^{1/\alpha} \sin{(\frac\pi\alpha)}.$
\end{Lemma}
\begin{proof}
We use the integral representation for the $r$th absolute moment \cite[Lemma 2]{Bahr+Esseen:1965}:
\begin{equation}    \label{eq:integral representation}
m_r \defeq \E[|W|^r] = \frac{\Gamma(r+1)}{\pi} \sin \Big(\frac{\pi r}{2}\Big)
\int_{\R} \frac{1-\Real\,\E [e^{\imag t W}]}{|t|^{r+1}} \, \dt.
\end{equation}
Set $K(r) \defeq \frac{\Gamma(r+1)}{\pi}\sin(\frac{r\pi}{2})$,
$B \defeq \Gamma(1-\alpha)\cos(\frac{\pi\alpha}{2})$ and
$C \defeq \Gamma(1-\alpha)\sin(\frac{\pi \alpha}{2})$.
Using Euler's identity $e^{\imag x}=\cos x +\imag\sin x$ in \eqref{eq:characteristic function}, we obtain
\begin{equation*}
\Real\,\E [e^{\imag tW}]
= \exp{(-B|t|^{\alpha})}\cos{(-C|t|^{\alpha}{\rm sgn}(t))}.
\end{equation*}
Substituting this into formula \eqref{eq:integral representation} yields
\begin{equation*}
m_r = 2 K(r) \int_0^{\infty} \frac{1-e^{-Bt^{\alpha}}\cos{(Ct^{\alpha})}}{t^{r+1}} \, \dt.
\end{equation*}
A change of variables ($u \defeq t^{\alpha}$) gives
\begin{align}
m_r
& =
\frac{2K(r)}{\alpha}\int_0^{\infty}\big(1-e^{-Bu}\cos{(Cu)}\big) u^{-1-r/\alpha} \, \du    \notag  \\
& =
\frac{2K(r)}{\alpha} \int_0^{\infty}\big(1-e^{-Bu}\big) u^{-1-r/\alpha} \, \du \notag  \\
&\hphantom{=}
+ \frac{2K(r)}{\alpha} \int_0^{\infty} e^{-Bu} \big(1-\cos{(Cu)}\big) u^{-1-r/\alpha} \, \du
 \eqdef
I_1+I_2.    \label{eq:I_1+I_2}
\end{align}
Integration by parts yields:
\begin{align*}
I_1
& =
\frac{2 K(r) B}{r}\int_0^{\infty}u^{-r/\alpha}e^{-Bu} \, \du    \\
& =
\frac{2 K(r) B^{r/\alpha}}{r} \Gamma\Big(1-\frac r \alpha\Big)
= -\frac{2 K(r) B^{r/\alpha}}{\alpha} \Gamma\Big(\!-\!\frac r\alpha\Big).
\end{align*}
According to \cite[Formula (3.945(2))]{Gradshteyn+Ryzhik:2000}, we have	%??? Matthias M. never checked this source! The rest of the proof is correct.
\begin{equation*}
I_2 = \frac{2 K(r)}{\alpha} \,
\Gamma\Big(\!-\!\frac r\alpha\Big)\Big(B^{r/\alpha}-|\Gamma(1-\alpha)|^{r/\alpha}
\cos\Big(\frac {\pi r}{2}-\frac{\pi r}{\alpha}\Big)\Big).
\end{equation*}
Now plugging in the values of $I_1$ and $I_2$ in \eqref{eq:I_1+I_2} gives
\begin{align*}
m_r
& =
-\frac{2 K(r)}{\alpha}\Gamma\Big(\!-\!\frac r\alpha\Big)|\Gamma(1-\alpha)|^{r/\alpha} \cos\Big(\frac {\pi r}{2}-\frac{\pi r}{\alpha}\Big) \\
& = \frac{2 K(r)}{r}\Gamma\Big(1\!-\!\frac r\alpha\Big)|\Gamma(1-\alpha)|^{r/\alpha} \cos\Big(\frac {\pi r}{2}-\frac{\pi r}{\alpha}\Big).
\end{align*}
The proof is complete.
\end{proof}

\begin{Lemma}		\label{Lem:power moments A4}
Let $W$ be a random variable with
\begin{equation*}
\E\big[e^{\theta W}\big] = E_{\alpha}\Big(\frac{\theta}{\Gamma(1-\alpha)}\Big)
\text{ for every } \theta \in \R
\end{equation*}
where $E_{\alpha}$ denotes the Mittag-Leffler function with parameter $\alpha \in (0,1)$.
Then, for any $r > 0$, we have
\begin{equation*}
\E [W^r]
= \frac{\Gamma(r+1)}{\Gamma(1-\alpha)^r \Gamma(r\alpha+1)}.
\end{equation*}	
\end{Lemma}
\begin{proof}
For $\alpha\in (0,1)$, let $S_\alpha$ denote a positive $\alpha$-stable random variable
with Laplace transform $\E[e^{-\lambda S_\alpha}]=\exp(-\Gamma(1-\alpha) \lambda^\alpha)$, $\lambda \geq 0$. We shall need the following integration formula for positive random variables $X$ and $s>0$
\begin{equation}	\label{eq:integration formula}
\E[X^{-s}] = \frac{1}{\Gamma(s)} \int_0^\infty t^{s-1} \E[e^{-tX}] \dt
\end{equation}
which follows from the fact that $\Prob(E/X > t) = \E[e^{-tX}]$ for all $t \geq 0$ where $E$ is a random variable with an exponential law with mean $1$ which is independent of $X$. Using \eqref{eq:integration formula} for $X = S_\alpha$ and $s=r\alpha$ gives
\begin{align*}
\E[S_\alpha ^{-r\alpha}] = \frac{1}{\Gamma(r\alpha)} \int_0^\infty t^{r\alpha-1} \E[e^{-tS_\alpha}] \dt
= \frac{1}{\Gamma(r\alpha)} \int_0^\infty t^{r\alpha-1} e^{-\Gamma(1-\alpha) t^\alpha} \dt
=  \frac{\Gamma(r+1)}{\Gamma(1-\alpha)^r \Gamma(r\alpha+1)}.
\end{align*}
This shows that the moment generating function of $S_\alpha^{-\alpha}$ is the same as that of $W$,
which proves that $S_\alpha^{-\alpha}$ has the same law as $W$.
In particular, $\E[W^r]=\E[S_\alpha ^{-r\alpha}]$ for all $r\geq 0$ which completes the proof.
\end{proof}

\section*{Acknowledgements}

\noindent
We would like to thank Yuan Li for pointing out an error in a previous version of this note.

\end{document}